\documentclass[12pt]{article}
\usepackage{amsmath, amssymb, eucal, latexsym}
\usepackage[cp1251]{inputenc}
\usepackage[english]{babel}
\usepackage{graphicx}
\def\Z{\Bbb Z}
\def\ve{\varepsilon}
\newtheorem{theorem}{Theorem}

\newtheorem{lemma}[theorem]{Lemma}

\newtheorem{note}[theorem]{Remark}
\newtheorem{definition}[theorem]{Definition}

\def\T{{\mathbb T}}
\def\oT{{\mathbf T}}
\def\C{{\mathbb C}}
\def\R{{\mathbb R}}
\def\Z{{\mathbb Z}}
\def\N{{\mathbb N}}
\def\l{\lambda}
\def\t{\tilde}
\def\v{\vec}
\def\_phi{\varphi}
\def\ov{\overline}
\def\L{\Lambda}

\begin{document}

\underline{.........................................................................................................................................................}    \\

\

 \begin{center}

                           \textbf{QUANTITATIVE VERSION OF BEURLING--HELSON THEOREM}                          \\

    \end{center}

\

 \begin{center}
                                                         S. V. KONYAGIN
\footnote{
The first author is supported by grant RFBR 14-01-00332
and grant Leading Scientific Schools N 3082.2014.1}                                                         ,
I. D. SHKREDOV\footnote{
The second author is supported by grant
mol\underline{ }a\underline{ }ved 12--01--33080.
}\\

    \end{center}

\

\begin{center}
    Abstract.
\end{center}

{\it \small
    It is proved that any continuous function $\_phi$ on the unit circle such that the sequence
    $\{ e^{in \_phi} \}_{n\in \Z}$ has small Wiener norm
    $\| e^{in \_phi} \|_{A} =  o \left( \frac{\log^{1/22} |n|}{(\log \log |n|)^{3/11}} \right)$,
    $|n| \to \infty$, is linear.
    Moreover, we get lower bounds for  Wiener norm of characteristic functions of subsets from $\Bbb Z_p$ in the case of prime $p$.}
\\

\section{Introduction}
\label{sec:introduction}

Let $A(\T)$ be the space of complex continuous functions $f$ on the unit circle
$\T=\R/(2\pi\Z)$, having absolutely  convergent series of its Fourier
coefficients. Equipping with the norm
$$
    \| f\|_{A(\T)} := \sum_{k\in \Z} |\hat f (k)|
$$
the space $A(\T)$ is a Banach algebra under the ordinary pointwise
multiplication of functions. Here we have used the notation
$$
    \hat f (k) = (2\pi)^{-1} \int_{\T} f(t) e^{-ikt} dt \,,\quad k \in \Z \,.
$$

Let $\_phi$ be a continuous map on $\T$ to itself, that is a continuous function $\_phi : \R \to \R$ such that
$\_phi (t+2\pi) \equiv \_phi(t) \pmod {2\pi}$.
A well--known Beurling--Helson theorem \cite{BH} asserts that if $\| e^{in\_phi} \|_{A(\T)} = O(1)$, $n\in \Z$
then the map $\_phi$ is linear. In other words, we have in this case
$\_phi(t) = \nu t + \_phi (0)$, $\nu \in \Z$.
The result gives a solution of a problem of P. L\'{e}vy on endomorphisms  of the algebra $A(\T)$ : any such an endomorphism is trivial, namely, it has the form $f(t) \to f(\nu t + t_0)$.
J.--P. Kahane \cite{Kahane,Kahane_book} conjectured that a weaker condition
$\| e^{in\_phi} \|_{A(\T)} = o(\log |n|)$, $|n| \to \infty$, implies that $\_phi$ is a linear function. Also he showed in \cite{Kahane_book} that any continuous piecewise linear but not linear map
 $\_phi : \T \to \T$ satisfies  $\| e^{in\_phi} \|_{A(\T)} \asymp \log |n|$, $|n| \to \infty$. Thus if the conjecture takes place then it is sharp.
 In paper \cite{Lebedev3} (see also \cite{Lebedev1,Lebedev2}) V.V. Lebedev obtained the first quantitative version of Beurling--Helson theorem.

\begin{theorem}
   Let $\_phi : \T \to \T$ be a continuous map.
   Suppose that
\begin{equation}\label{cond:Lebedev}
    \| e^{in \_phi} \|_{A(\T)} = o \left( \left( \frac{\log \log |n|}{(\log \log \log |n|)} \right)^{1/12} \right) \,, \quad n\in \Z,\, \quad |n| \to \infty\,.
\end{equation}
    Then for some $\nu \in \Z$ the following holds $\_phi (t) = \nu t+\_phi(0)$.
\label{t:Lebedev}
\end{theorem}

Let us formulate the main result of the paper.

\begin{theorem}
   Let $\_phi : \T \to \T$ be a continuous map.
   Suppose that
\begin{equation}\label{cond:e_phi_intr}
    \| e^{in \_phi} \|_{A(\T)} = o \left( \frac{\log^{1/22} |n|}{(\log \log |n|)^{3/11}} \right) \,, \quad n\in \Z,\, \quad |n| \to \infty\,.
\end{equation}
    Then for some $\nu \in \Z$ the following holds $\_phi (t) = \nu t+\_phi(0)$.
\label{t:main_intr}
\end{theorem}

In our proof we develop the approach from \cite{Lebedev1,Lebedev2,Lebedev3}. The main advantage is connected with the notion of the additive dimension (\S2). After discretization of the problem an appropriate upper bound for the dimension  allows us to
consider the values of the function $\varphi$ not in all points of the
lattice but on its small subset.
Besides we win a little using the results of T.~Sanders \cite{Sanders} on lower bounds of Wiener norm
of characteristic functions of large subsets from $\Bbb Z_p$ (\S3) in the case of prime $p$.
As a byproduct we obtain the best possible lower estimates for Wiener norm of characteristic functions of small subsets of $\Bbb Z_p$.

\section{Dissociated sets and its generalizations}
\label{sec:dissoc}

Let $G$ be a compact abelian group and $\hat G$ be the dual group, that is the discrete abelian group of continuous homomorphisms $\gamma:\,G\to S^1$, where $S^1=\{z\in\Bbb C:\,|z|=1\}$. It is a well--known fact that there is the Haar measure  $\mu_G$  on
$G $.
We suppose that the measure $\mu_G$ is normalized such as $\mu_G(G)=1$. Denote
the Fourier transform of a function $f\in L^1(G)$, that is an arbitrary integrable complex function on $G$ to be a new function
$\hat f\in l^\infty(\hat G)$
defined by the formula
$$\hat f(\gamma)=\int_{x\in G}f(x)\overline{\gamma(x)}d\mu_G(X).$$
It is a well--known fact that for any functions $f,g\in L^1(G)$ the following holds
\begin{equation}\label{scal_prod}
\int_{x\in G}f(x)\overline g(x)d\mu_G(X)=\sum_{\gamma\in \hat G}\hat f(\gamma)\overline{\hat g(\gamma)}.
\end{equation}

Let
$$A(G)=\{f\in L^1(G):\,\| \hat f \|_1<\infty\} \,,$$
where $\| \hat f \|_1=\sum_\gamma |\hat f(\gamma)|$,
and define a norm on $A(G)$ as $\|f\|_{A(G)}=\|\hat f\|_1$. Note that any function belonging to $A(G)$ is a continuous function on $G$.

By $\chi_S$, $S\subset G$ we denote the characteristic function of some set $S$.

Our main object is the group $\Z_p=\Z/p\Z$, where $p$ is a prime number.

Given an arbitrary function $f:G\to \C$ denote the quantity
$$
    \oT_k (f) := \sum_{x_1+\dots+x_k = x'_1+\dots+x'_k} f(x_1) \dots f(x_k) \ov{f(x'_1)} \dots \ov{f(x'_k)} \,.
$$
In the case of finite group $G$, we have
$$
    \oT_k (f) =
        |G|^{2k-1} \sum_\gamma |\hat {f} (\gamma)|^{2k} \,.
$$

The main idea of the proof of our Theorem \ref{t:main_intr} can be demonstrated by the following arguments which connects Wiener norm of the characteristic function and its additive dimension.

 Recall that a set $\Lambda = \{ \lambda_1, \dots, \lambda_{|\Lambda|} \}$ is called {\it dissociated} if any equality
\begin{equation}\label{f:diss_classical}
    \sum_{\l \in \Lambda}\ve_{\l}\l=0 \,,
\end{equation}
where $\ve_{\l}\in \{-1,0,1\}$ implies that all  $\ve_{\l}$ are equal to zero.
It is well--known (see e.g. \cite{Rudin_book}) that $\oT_k (\Lambda) \le (Ck)^k |\Lambda|^k$, $k\ge 2$,
where $C>0$ is an absolute constant.
The size of a maximal dissociated subset of a set $S\subseteq G$ is called
{\it the additive dimension} of $S$ and is denoted as $\mathrm{dim} (S)$.

\begin{theorem}
    Let $G$ be a finite abelian group,
    $S\subseteq G$ be any set, $\| \chi_S \|_{A(G)} \le K$,
    and
    $K^2 \le |S|$.
    Then $\mathrm{dim} (S) \ll K^2 (1+\log (|S|/K^2))$.
\label{p:dim_A(G)}
\end{theorem}

To get Theorem \ref{p:dim_A(G)} we need the following lemma.

\begin{lemma}\label{T_k_est}
Let $G$ be a finite abelian group, $f$ be a complex function on $G$,
$\| f \|_{A(G)} \le K$, $Q\subset S\subset G$, $|f(x)|=1$ for $x\in S$, and $f(x)=0$ for $x\in G\setminus S$,
$g(x)=\chi_Q(x)f(x)$, $k\in\N$. Then
$$\oT_k (g) \ge \frac{|Q|^{2k}}{|S| K^{2k-2}}.$$
\end{lemma}

{\bf Proof.}
Using identity~(\ref{scal_prod}), H\"{o}lder inequality and the assumption $\| S \|_{A(G)} \le K$, we have
$$
    |Q|^{2k} = \left( \sum_x g(x) \overline{f(x)} \right)^{2k}
            \le \left( |G| \sum_\gamma |\hat g (\gamma)| |\hat f (\gamma)| \right)^{2k}
        \le
$$
$$
        \le
            |G|^{2k} \sum_\gamma |\hat g (\gamma)|^{2k}
\left( \sum_\gamma |\hat f (\gamma)|^{2k/(2k-1)} \right)^{2k-1}\le
$$
$$
                \le
                    |G| \oT_k (g) \left( \sum_\gamma |\hat f (\gamma)| \right)^{2k-2}
                        \cdot \sum_\gamma |\hat f (\gamma)|^2
                        \le
                    \oT_k (g) K^{2k-2} |S| \,.
$$
This completes the proof.
$\hfill\Box$

\bigskip

{\bf Proof of Theorem~\ref{p:dim_A(G)}.} Put $f=\chi_S$.
Letting $Q$ be a maximal dissociated set $\Lambda \subseteq S$, we get from Lemma~\ref{T_k_est}
and Rudin's inequality
$$
    (Ck)^k |\Lambda|^k \ge \frac{|\Lambda|^{2k}}{|S| K^{2k-2}} \,,
$$
where $C>0$ is an absolute constant.
Taking $k=2+[\log (|S|/K^2)]$, we obtain $|\Lambda| \ll  K^2 (1+\log (|S|/K^2))$.
This completes the proof.
$\hfill\Box$

\bigskip

\begin{note}
\label{Knote}
    If $S\subseteq G$ is any set then the Parseval identity implies that
    $\| S \|_{A(G)} \le |S|^{1/2}$. Thus the condition $K^2 \le |S|$ of Theorem \ref{p:dim_A(G)} is not burdensome.
    Note also that the first multiple $K^2$ in the estimate of dimension of $S$ in Theorem~\ref{p:dim_A(G)} cannot be replaced by $K^{2-\varepsilon}$, $\varepsilon>0$, generally speaking. Indeed, it is sufficient to consider a random set $S$, for example.
\end{note}

\bigskip

In paper \cite{Shkr_large_exp} another version of the definition
of a dissociated set different from the classical one (\ref{f:diss_classical}) was considered.

\begin{definition}
   Let $k,s$ be positive integers.
   A set $\L = \{ \l_1, \dots, \l_{|\L|} \}\subset\Z_N$
   belongs to the family  $\L(k,s)$ if any equality
   \begin{equation}\label{}
        \sum_{i=1}^{|\L|} s_i \l_i = 0
            \,, \quad \l_i \in \L \,, \quad s_i \in \Z \,, \quad
        |s_i| \le s\,, \quad \sum_{i=1}^{|\L|} |s_i| \le 2k \,,
   \end{equation}
   implies that all $s_i$ are equal to zero.
\label{def:diss_k,s}
\end{definition}

In the same paper the following results was proved (see \cite{Shkr_large_exp}, Statement 1).

\begin{lemma}
            Let $k$, $s$ be positive integers, $\L$ be a set from the family $\L(k,s)$,
            and $|\L| \ge k$. Then
        \begin{equation}\label{f:T_k_up_estimate}
            \oT_k (\L) \le 2^{3k} k^k |\L|^k \max \left\{ 1, \left( \frac{k}{|\L|} \right)^k |\L|^{k/s} \right\} \,.
        \end{equation}
\label{l:Lambda_k}
\end{lemma}

Using Definition \ref{def:diss_k,s} as well as Lemma \ref{l:Lambda_k} one can obtain an analog of Theorem \ref{p:dim_A(G)}, where the dimension of a set is defined as the size of its a maximal subset from the family $\Lambda(k,s)$. We do not need in this sharper result.

Given $u\in \T$ put
$$\|u\|=\frac1{2\pi}\inf_{v\in\R, v\equiv u(\bmod 2\pi)}|v|.$$
We need in an analog of Definition \ref{def:diss_k,s} relatively to a function
$\_phi  : \T \to \T$.
For any $x\in \Z_N$ we write $\_phi^* (x)$ for $\_phi (x/N)$.

\begin{definition}
   Let $k,s$ be positive integers, $\eta \in (0,1]$ be a real number, and
   $\_phi : \T \to \T$ be any function.
   A set $\L = \{ \l_1, \dots, \l_{|\L|} \}\subset\Z_N$ belongs to the family
   $\L^{\_phi, \eta} (k,s)$ if any equality
   \begin{equation}\label{}
        \sum_{i=1}^{|\L|} s_i \l_i = 0
            \,, \quad \l_i \in \L \,, \quad s_i \in \Z \,, \quad
        |s_i| \le s\,, \quad \sum_{i=1}^{|\L|} |s_i| \le 2k \,,
   \end{equation}
   and inequality
   \begin{equation}\label{}
        \| \sum_{i=1}^{|\L|} s_i \_phi^* (\l_i) \| \le \eta
   \end{equation}
   imply that all $s_i$ are equal to zero.
\label{def:diss_k,s_phi}
\end{definition}

Clearly, any subset of a set from the family $\L^{\_phi, \eta} (k,s)$ automatically belongs to the family.

 Having a function $\_phi : \T \to \T$, a real number $\eta \in (0,1]$ and an arbitrary set $S \subseteq \Z_N$ let us define the quantity $\oT_k^{\_phi,\eta} (S)$ as
$$
    \oT_k^{\_phi,\eta} (S) :=
        | \{ (x_1,\dots, x_{2k}) \in S^{2k}  ~:~ x_1 + \dots + x_k = x_{k+1} + \dots + x_{2k} \,,
$$
$$
        \| \_phi^* (x_1) + \dots + \_phi^* (x_k) - \_phi^* (x_{k+1}) - \dots - \_phi^* (x_{2k}) \| \le \eta \} | \,.
$$

In a similar way we obtain a statement on an
upper bound
of quantity  $\oT_k^{\_phi,\eta} (\L)$ for dissociated sets $\L$ from the family
$\L^{\_phi, \eta} (k,s)$.

\begin{lemma}
            Let $k$, $s$ be positive integers, $s\ge 5$,
        $\L$ be a set from the family $\L^{\_phi,\eta} (2k,s)$.
        Then
        \begin{equation}\label{f:T_k_up_estimate_phi}
            \oT^{\_phi,\eta/2}_k (\L) \le 2^{4k+2} k^{k+1} |\L|^k \max \left\{ 1, \left( \frac{k}{|\L|} \right)^k |\L|^{4k/s} \right\} \,.
        \end{equation}
\label{l:Lambda_k_phi}
\end{lemma}

The proof is close to the arguments  from  \cite{Shkr_large_exp},
the only difference is that we estimate quantity $\oT^{\_phi}_k$ directly and, hence, do not use the mean value of multidimensional convolution (more precisely see {\cite{Shkr_large_exp}).
The approach gives slightly weaker bounds but our implications to Kahane's problem do not depend on this.

{\bf Proof of Lemma \ref{l:Lambda_k_phi}}.
Let $\sigma = \oT^{\_phi,\eta/2}_k (\L)$ be the number of tuples $(\l_1,\dots,\l_{2k}) \in \L^{2k}$
satisfying the equation
\begin{equation}\label{f:eq_phi}
    \l_1+\dots+\l_k = \l_{k+1} + \dots + \l_{2k}
\end{equation}
as well as the inequality
   \begin{equation}\label{f:ineq_phi}
        \| \_phi^* (\l_1) + \dots + \_phi^* (\l_k) - \_phi^* (\l_{k+1}) - \dots - \_phi^* (\l_{2k}) \| \le \eta/2 \,.
   \end{equation}
Reducing equal terms in the equation, rewrite formula (\ref{f:eq_phi}) as
\begin{equation}\label{tmp:04.12.2013_1}
    \sum_{j=1}^l s_j \t{\l}_j = 0 \,, \quad s_j \in \Z \setminus \{ 0 \} \,,
\end{equation}
where all $\t{\l}_j\in \L$ are different and $\sum_{j=1}^l |s_j| = 2t \le 2k$.
Split sum (\ref{tmp:04.12.2013_1}) onto two sums
\begin{equation}\label{tmp:04.12.2013_2}
    0 = \sum_{j=1}^l s_j \t{\l}_j = \sum_{j\in G} s_j \t{\l}_j + \sum_{j\in B} s_j \t{\l}_j \,,
\end{equation}
where $G=\{ j ~:~ |s_j| \le s/2 \}$ and $B=\{ j ~:~ |s_j| > s/2 \}$.
Clearly, the sets $G$ and $B$ depend on the sequence $\v{s} = (s_1,\dots,s_l)$.

Fix $l,t$ and a sequence $\v{s} = (s_1,\dots,s_l)$ and estimate the number of solutions of equality (\ref{tmp:04.12.2013_2}).
Consider the solutions of (\ref{tmp:04.12.2013_2}) with fixed $\t{\l}_j$, $j\in B$.
Take any two of such solutions $(\t{\l}'_1, \dots, \t{\l}'_l)$, $(\t{\l}''_1, \dots, \t{\l}''_l)$.
We have
\begin{equation}\label{tmp:04.12.2013_3}
    \sum_{j\in G} s_j \t{\l}'_j - \sum_{j\in G} s_j \t{\l}''_j = 0
\end{equation}
and the triangle inequality implies that
$$
    \| \sum_{j\in G} s_j \_phi^* (\t{\l}'_j) - \sum_{j\in G} s_j \_phi^* (\t{\l}''_j) \| \le \eta \,.
$$
By the definition of the set $G$, we get that all $s_j$ in the formula above does not exceed $s$.
Hence by the definition of the family $\L^{\_phi,\eta} (2k,s)$
the tuple $\{ \t{\l}''_j \}_{j\in G}$ is a permutation of the tuple $\{ \t{\l}'_j \}_{j\in G}$.
The number of such permutations of a fixed tuple equals  $\frac{(2t)!}{|s_1|! \dots |s_l|!}$.
Whence the number of the solutions of equation (\ref{tmp:04.12.2013_2}) with fixed $l,t$ and $\v{s}$
does not exceed
    \begin{equation}\label{tmp:06.01.2014_1}
        \frac{(2t)!}{|s_1|! \dots |s_l|!} |\L|^{|B (\v{s})|} \,.
    \end{equation}

For any $t$, $0\le t \le k$ put $r_t=2t /([s/2]+1)$.
Clearly, for any tuple $\v{s} = (s_1,\dots,s_l)$, $\sum_{j=1}^l |s_j| = 2t$ the following holds $|B(\v{s})| \le r_t$.
Using estimate (\ref{tmp:06.01.2014_1}) and inequality
$$
  \sum_{l=1}^{2t} \sum_{s_1+\dots+s_l = 2t,\, s_j> 0}\, \frac{(2t)!}{s_1! \dots s_l!} \le
  \sum_{s_1+\dots+s_{2t} = 2t,\, s_j\ge 0}\, \frac{(2t)!}{s_1! \dots s_l!} = (2t)^{2t} \,,
$$
we obtain the required upper bound for the quantity $\sigma$
$$
    \sigma
        \le
            \sum_{t=0}^k \sum_{b=0}^{r_t}\, \sum_{\v{s} ~:~ |B(\v{s})| = b} \frac{(2t)!}{|s_1|! \dots |s_l|!} |\L|^b
            (k-t)! |\L|^{k-t}
                \le
                    2^{4k} |\L|^k \sum_{t=0}^k \sum_{b=0}^{r_t}\, t^{2t} k^{k-t} |\L|^{b-t}
                        \le
$$
$$
    \le
        2^{4k+1} k^k |\L|^k \sum_{t=0}^k \left( \frac{k}{|\L|^{1-4/s}} \right)^t
            \le
                2^{4k+2} k^{k+1} |\L|^k \max \left\{ 1, \left( \frac{k}{|\L|} \right)^k |\L|^{4k/s} \right\} \,.
$$

This concludes the proof of Lemma \ref{l:Lambda_k_phi}.
$\hfill\Box$

\section{On Fourier transform of characteristic functions of subsets of $\Bbb Z_p$}
\label{sec:Sanders}

As in the previous section $G$ is a compact abelian group and $\hat G$ is the dual group.

For a measurable set $E\subset G$ of positive measure, a continuous function $f$ on $G$,
and a real number $1\le q<\infty$ we define
$$\|f\|_{L^q(E)}=\left((\mu_G(E))^{-1}\int_{x\in E}|f(x)|^qd\mu_G(x)\right)^{1/q},$$
and
$$\|f\|_{L^\infty(E)}=\sup_{x\in E}|f(x)|.$$

In  \cite{GK} the following question was studied.
Let $A$ be a subset of $\Z_p$. How to obtain a lower bound for the quantity  $\|\chi_A\|_{A(\Z_p)}$ in terms of $|A|$?
Note that because of
$$\|\chi_{\Z_p\setminus A}\|_{A(\Z_p)}=\|\chi_A\|_{A(\Z_p)}+(1-2|A|/p)$$
it is sufficient to consider the case
$$|A|<p/2.$$
From the results of \cite{GK} it follows that in the situation case one has
$$\|\chi_A\|_{A(\Z_p)}\gg\frac{|A|}{p}\left(\frac{\log p}{\log\log p}\right)^{1/3}.$$
This estimate can be improved using a theorem of T.~Sanders \cite{Sanders}, which we describe below.

Given $z\in S^1$ put
$$\|z\|=\frac1{2\pi}\inf_{z\in\Z}|2\pi n+\arg z|.$$
Thus $\|z\|$ is small if $z$ is close to 1. For a finite nonempty set $\Gamma \subseteq \hat G$ and $\delta\in(0,1]$ define the Bohr set
$$B(\Gamma,\delta)=\{x\in G:\,\|\gamma(x)\|\le\delta\quad\text{for all}\quad \gamma\in\Gamma\}.$$
Shifts of Bohr sets are called Bohr neighborhoods. We write $d=|\Gamma|$.

In the theory of Bohr sets the following result plays an important role (see e.g. \cite{Sanders}, Lemma~6.2).
\begin{lemma}
\label{sizeBohrset}
One has $\mu_G(B(\Gamma,\delta))\ge\delta^d$.
\end{lemma}

For fixed $\Gamma$ and $\delta$ we denote by $\beta (x)$ the function on $G$ such that
$\beta (x)$  equals $1/\mu_G(B(\Gamma,\delta))$ for  $x\in B(\Gamma,\delta)$ and
zero otherwise.
Given two functions $f,g\in L^1(G)$ define its convolution as
$$f*g(x)=\int_{y\in G} f(y)g(x-y)d\mu_G(y).$$
Note that the convolution of a continuous function and an integrable function is also continuous.
It is well--known that
\begin{equation}
\label{four_conv}
\widehat{f*g}(\gamma)=\hat f(\gamma)\hat g(\gamma).
\end{equation}

\bigskip

The following theorem is a main result of paper \cite{Sanders}.
\begin{theorem}
\label{Sanders}
Let $G$ be a compact abelian group, $f\in A(G)$, $f\not\equiv0$ and $\ve\in(0,1]$. Put
$A_f=\|f\|_{A(G)}\|f\|_\infty^{-1}$. Then there is  a Bohr set $B(\Gamma,\delta)$ such that
$$d\ll \ve^{-2}A_f(1+\log A_f)(1+\log(\ve^{-1}A_f))$$
and
$$\log(\delta^{-1})\ll\ve^{-2}A_f(1+\log(\ve^{-1}A_f)),$$
and also a smaller Bohr set $B(\Gamma,\delta')$, $\delta'\gg\ve\delta/d$ with
$$\sup_{x\in G}\|f*\beta-f*\beta(x)\|_{L^\infty(x+B(\Gamma,\delta'))}\le\ve\|f\|_{L^\infty(G)}$$
and
$$\sup_{x\in G}\|f-f*\beta\|_{L^2(x+B(\Gamma,\delta'))}\le\ve\|f\|_{L^\infty(G)}.$$
\end{theorem}

Using Theorem~\ref{Sanders} we obtain a result.
\begin{theorem}
\label{Char}
Let $p$ be a prime number, $A\subset\Z_p$, $0<\eta=|A|/p<1/2$. If $\eta\ge(\log p)^{-1/4}(\log\log p)^{1/2}$ then
$$\|\chi_A\|_{A(\Z_p)}\gg(\log p)^{1/2}(\log\log p)^{-1}
\eta^{3/2} \left(1 + \log\left(\eta^2(\log p)^{1/2}(\log\log p)^{-1}\right)\right)^{-1/2},$$
and if  $\eta<(\log p)^{-1/4}(\log\log p)^{1/2}$ then
$$\|\chi_A\|_{A(\Z_p)}\gg\eta^{1/2}(\log p)^{1/4}(\log\log p)^{-1/2}.$$
\end{theorem}

{\bf Proof of Theorem~\ref{Char}}. Putting  $G=\hat G=\Z_p$ and
$$f_0=\chi_A-\eta \,,$$
we have
\begin{equation}
\label{f0average}
\hat f_0(0)=\int_G f_0d\mu_G=0.
\end{equation}
Further let $k=[1/(2\eta)]$ and consider a random function
$$f(x)=\sum_{j=1}^k f_0(x-x_j) \,,$$
where $x_1,\dots,x_k$ are uniformly distributed independent variables on $G$.
For $\gamma\in\hat G$, we have
$$\hat f(\gamma)=\hat f_0 (\gamma)\sum_{j=1}^k\gamma(x_j).$$
Because of~(\ref{f0average})
\begin{equation}
\label{faverage}
\hat f(0)=0.
\end{equation}
In the case $\gamma\neq0$ the second moment of  $\hat f(\gamma)$ is
$${\bold E}|\hat f(\gamma)|^2 = k|\hat f_0(\gamma)|^2.$$
Thus by the Cauchy--Schwartz inequality, we get
$${\bold E}|\hat f(\gamma)| \le\sqrt k|\hat f_0(\gamma)|.$$
Hence
$${\bold E}\|f\|_{A(G)} \le\sqrt k\|f_0\|_{A(G)}.$$
Now fix $x_1,\dots,x_k$ such that
\begin{equation}
\label{f_via_f_0}
\|f\|_{A(G)} \le\sqrt k\|f_0\|_{A(G)}.
\end{equation}
We will use two trivial estimates for $\|f\|_\infty$:
\begin{equation}
\label{infty_via_trian}
\|f\|_\infty\le k
\end{equation}
and
\begin{equation}
\label{infty_via_Anorm}
\|f\|_\infty\le\|f\|_{A(G)}.
\end{equation}
If the condition
\begin{equation}
\label{cond_nontr}
(\delta')^d>1/p
\end{equation}
holds then by Lemma~\ref{sizeBohrset}, one has $\mu_G(B(\Gamma,\delta'))>1/p$, that is
the set $B(\Gamma,\delta')$ contains a nonzero element $x_0$. In view of~(\ref{four_conv})
and~(\ref{faverage})
$$\sum_{j=0}^{p-1} (f*\beta) (jx_0)=p (f*\beta)\hat{}\, (0)=p\hat f(0) \cdot \hat\beta(0)=0.$$
Hence there is an element $x=jx_0$ such that  $(f*\beta) (x) (f*\beta) (x+x_0)\le0$. Further by Theorem~\ref{Sanders}, we obtain  $|f*\beta(x+x_0)-f*\beta(x)|\le\ve\|f\|_{L^\infty(G)}$. It follows that
$|f*\beta(x_1)|\le0.5\ve\|f\|_{L^\infty(G)}$ for either $x_1=x$ or $x_1=x+x_0$. Whence for any  $x\in x_1+B(\Gamma,\delta')$,
we have $|f*\beta(x)|\le1.5\ve\|f\|_{L^\infty(G)}$. Using Theorem~\ref{Sanders} once more time, we get
$$\|f-f*\beta\|_{L^2(x_1+B(\Gamma,\delta'))}\le\ve\|f\|_{L^\infty(G)}.$$
Thus there is an element $x\in x_1+B(\Gamma,\delta')$ with
$|f(x)-f*\beta(x)|\le\ve\|f\|_{L^\infty(G)}$. Finally,
\begin{equation}
\label{est_elem}
|f(x)|\le2.5\ve\|f\|_{L^\infty(G)}.
\end{equation}
On the other hand, by the definition of the number $k$  and the function $f$ it follows that $|f(x)|>0.25$ for any
$x\in G$. Hence if inequality
\begin{equation}
\label{ineq_eps}
\ve\|f\|_{L^\infty(G)}\le0.1
\end{equation}
takes place then we get a contradiction.

We use the notation
$$u=\|f\|_{A(G)},\quad v=\|f\|_\infty.$$
Note that $u\ge v$. Put
$$\ve=0.1v^{-1}.$$
Thus inequality~(\ref{ineq_eps}) holds. Consider two cases.

Case 1: $\eta\ge(\log p)^{-1/4}(\log\log p)^{1/2}$. Let
$$u_1=\eta (\log p)^{1/2}(\log\log p)^{-1}
\left(1+\log\left(\eta^2(\log p)^{1/2}(\log\log p)^{-1}\right)\right)^{-1/2},$$
$$u_0=cu_1,$$
where the constant $c\in(0,1)$ depends on the constants in the signs $\ll$ in Theorem~\ref{Sanders}.
We will use the fact that our assumption on $\eta$ implies that the logarithm in the second multiple of quantity $u_1$ is nonnegative and, hence, the multiple does not exceed one.
Suppose that
\begin{equation}
\label{supp1_u}
u\le u_0.
\end{equation}
Our aim is to check inequality~(\ref{cond_nontr}) which gives us a contradiction with (\ref{supp1_u}).

Putting
\begin{equation}
\label{def_v0}
v_0=\min(k,u_0),
\end{equation}
we get in view of~(\ref{infty_via_trian}) that
\begin{equation}
\label{supp1_v}
v\le v_0.
\end{equation}
Note that
$$\ve^{-2}A_f\ll v^2(u/v)=uv\le u_0v,$$
$$1+\log A_f\ll 1+\log(u_0/v),$$
$$1+\log(\ve^{-1}A_f)\ll 1+\log u_0 \ll\log\log p.$$
Let us estimate $d$.
Multiplying the last three bounds, we obtain
$$d\ll u_0(\log\log p) v(1+\log(u_0/v)).$$
Because of the function $v\to v(1+\log(u_0/v))$ is increasing for $v\in(0,u_0)$, we have by
(\ref{supp1_v}) that
\begin{equation}
\label{est1_d}
d\ll u_0v_0(1+\log(u_0/v_0))\log\log p.
\end{equation}
Further
$$\log(\delta^{-1})\ll u_0v_0\log\log p.$$
It is easy to see that $\log d\ll\log\log p$.
Hence
\begin{equation}
\label{est1_logdelta'}
\log((\delta')^{-1})\le \log (\delta^{-1})+\log d+\log(\ve^{-1})\ll u_0v_0\log\log p.
\end{equation}
Combining (\ref{est1_d}) and (\ref{est1_logdelta'}), we obtain
\begin{equation}
\label{est1_prod}
d\log((\delta')^{-1})\ll u_0^2v_0^2(1+\log(u_0/v_0))(\log\log p)^2.
\end{equation}

If $v_0=k$ then
$$d\log((\delta')^{-1})\ll c^2u_1^2k^2(1+\log(u_1/k))(\log\log p)^2.$$
(The constants in the sign $\ll$ depend on the corresponding constants in Theorem~\ref{Sanders},
and the choice of the constant $c$ is in our hands). Thus, to prove (\ref{cond_nontr}) it is sufficient to check that
\begin{equation}
\label{comp_logp}
u_1^2k^2(1+\log(u_1/k))(\log\log p)^2\ll\log p.
\end{equation}
By the definition of the quantity $u_1$, we get
\begin{equation}
\label{estgamma}
u_1\le \eta (\log p)^{1/2}(\log\log p)^{-1}.
\end{equation}
Hence, because of $k=[1/(2\eta)]$, we have
$$1+\log(u_1/k)\ll 1+\log\left(\eta^2(\log p)^{1/2}(\log\log p)^{-1}\right).$$
The last bound and the definition of $u_1$ imply that
$$u_1^2k^2(1+\log(u_1/k))(\log\log p)^2\ll\log p $$
and (\ref{comp_logp}) is proved.

In the situation $v_0=u_0<k$ inequality~(\ref{est1_prod}) gives us
\begin{equation}
\label{est2_prod}
d\log((\delta')^{-1})\ll u_0^4(\log\log p)^2,
\end{equation}
hence
$$d\log((\delta')^{-1})\ll c^2u_1^2k^2(\log\log p)^2.$$
In view of~(\ref{estgamma})
$$u_1^2k^2(\log\log p)^2\ll\log p,$$
and estimate~(\ref{cond_nontr}) holds again.

We get a contradiction with~(\ref{supp1_u}). Thus, $u>u_0$.
Using inequality~(\ref{f_via_f_0}), we obtain
\begin{equation}
\label{estviau0}
\|f_0\|_{A(G)}\ge u_0/\sqrt k,
\end{equation}
and the required result follows.

Case 2: $\eta<(\log p)^{-1/4}(\log\log p)^{1/2}$. Let
$$u_1=(\log p)^{1/4}(\log\log p)^{-1/2},\quad u_0=cu_1.$$
By our choice of the parameter $\eta$, we have $k>u_0$ and hence $v_0 = u_0$.
It is easy to see that the new choice of parameters preserves all calculations
in lines (\ref{supp1_v})---(\ref{comp_logp}).
Because of $v_0 = u_0$, using inequality~(\ref{est2_prod}), we obtain (\ref{cond_nontr}).

Thus we get a contradiction with~(\ref{supp1_u}) again. It follows that $u>u_0$.
Applying inequality~(\ref{estviau0}), we have the required result at the second case.
This completes the proof.
$\hfill\Box$

Theorem \ref{Char} is nontrivial if our subset $A$ is large, that is
$$|A|p^{-1}(\log p)^{1/2} (\log\log p)^{-1} \to\infty $$
(and of course $|A|<p/2$). One can hope that for any $2\le|A|<p/2$
the following holds
\begin{equation}
\label{Litt_conj}
\|\chi_A\|_{A(\Z_p)}\gg\log |A|.
\end{equation}
It is easy to see that the bound $\log|A|$ is attained in the case $A=\{1,\dots,|A|\}$.
For sets $B\subset\Z$ an analog of (\ref{Litt_conj}) is a well--known fact,
namely, it is proved in \cite{Kon}  and \cite{MGPS} that if $B\subset\Z$, $2\le|B|<\infty$ then
\begin{equation}
\label{Litt_thm}
\int_{-\pi}^\pi \left|\sum_{b\in B} e^{ibx}\right|dx\gg\log|B|.
\end{equation}
The following statement which is a consequence of Theorem~\ref{p:dim_A(G)} gives us an optimal lower bound for Wiener norm of small subsets $A$.

\begin{theorem}
\label{Charsmall}
Let $p$ be a prime number, $A\subset\Z_p$, and
$$2\le|A|\le\exp\left((\log p/\log\log p)^{1/3}\right).$$
Then~(\ref{Litt_conj}) takes place.
\end{theorem}

{\bf Proof of Theorem~\ref{Charsmall}.}
One can suppose that $p$ is a sufficiently large prime number.
We identify elements of $\Z_p$ with corresponding integers having the least absolute values.

Suppose that
\begin{equation}
\label{K_suppos}
K:=\|\chi_A\|_{A(\Z_p)}\le c(\log p/\log\log p)^{1/3} \,,
\end{equation}
where $c$ is an appropriate positive constant.
Note that if inequality (\ref{K_suppos}) does not hold then the required result
follows immediately from the condition on the size of $A$.

Applying Theorem~\ref{p:dim_A(G)}, we have in view of Remark~\ref{Knote} that
$$d:=\mathrm{dim} (A)\le\log p/\log\log p \,,$$
provided by the constant $c$ is sufficiently small.
Let $\Lambda$ be a maximal dissociated subset of $A$. Thus $|\Lambda|=d$. By Dirichlet theorem
there is a positive integer $q<p$ such that for any $\l\in\Lambda$ one has
$$|q\l|\le p^{1-1/d}\le p/\log p \,.$$
(Recall that $|q\l|$ is the absolute value of the correspondent integer from $[-p/2 , p/2]$). Put
$$B=\{qa:\,a\in A\}.$$
Clearly,
\begin{equation}
\label{A_as_B}
\|\chi_A\|_{A(\Z_p)}=\|\chi_B\|_{A(\Z_p)}.
\end{equation}
By our choice of the set $\Lambda$ any element $a\in A$ can be represented as
$$a=\sum_{\l \in \Lambda}\ve_{\l}\l \,,$$
where $\ve_{\l}\in \{-1,0,1\}$. Hence an arbitrary $b\in B$ is represented as
$$b=\sum_{\l \in \Lambda}\ve_{\l}q\l.$$
It follows that
$$|b|\le dp/\log p\le p/\log\log p.$$
Whence $B\subseteq [-p/3,p/3]$. Applying Theorem 7.28, chapter 10, book \cite{Zyg},  we get
\begin{equation}
\label{Zyg}
\|\chi_B\|_{A(\Z_p)}\gg\int_{-\pi}^\pi \left|\sum_{b\in B} e^{ibx}\right|dx.
\end{equation}
Combining formula~(\ref{A_as_B}) and inequalities~(\ref{Zyg}), (\ref{Litt_thm}),
we obtain
$$\|\chi_A\|_{A(\Z_p)}\gg\log|B|=\log|A|,$$
as required.
$\hfill\Box$

\bigskip

In a future paper we are going to consider the question on lower bounds of $\|\chi_A\|_{A(\Z_p)}$ in terms of the cardinality of the set
$A$
in the case of subsets $A\subset\Z_p$ having medium size.

\section{Proof of Theorem \ref{t:main_intr}}
\label{sec:Lebedev}

It is obvious that we can consider just nonnegative $n$ in formula (\ref{cond:e_phi_intr}).
Put
$$
    \Theta (n) = \max_{l=0,1,\dots,n} \| e^{il\_phi} \|_{A(\T)} \,.
$$
Clearly, $\Theta (n)$, $n=0,1,2,\dots$ is a nondecreasing sequence and  $\Theta (0) = 1$.
By the assumption of Theorem \ref{t:main_intr}, we have
\begin{equation}\label{cond:main_22_1}
    \Theta (n) = o \left( \frac{\log^{1/22} n}{(\log \log n)^{3/11}} \right) \,, \quad n \to \infty \,.
\end{equation}

Combining Lebedev's arguments from \cite{Lebedev3} and Theorem \ref{Char}, we obtain the following auxiliary statement.

\begin{lemma}
    Let $\_phi : \T \to \T$ be a continuous map.
    Let also $Q$ be a positive integer, $N$ be a prime number, $N\le Q$ and for all $x\in \Z_N$ the following holds
\begin{equation}\label{cond:Lebedev_extract}
    \| Q \_phi^* (x) \| \le 1/N \,.
\end{equation}
    Put $\eta = 1/(64 \Theta^2 (Q))$ and let $\eta\ge(\log N)^{-1/4}(\log\log N)^{1/2}$.
    Then
$$
    \Theta^2 (Q)
        \gg (\log N)^{1/2}(\log\log N)^{-1} \eta^{3/2} \left(1 + \log\left(\eta^2(\log N)^{1/2}(\log\log N)^{-1}\right)\right)^{-1/2} \,.
$$
\label{l:Lebedev_extract}
\end{lemma}

Now let us prove our main Theorem \ref{t:main_intr}.
Again we use the arguments of Lebedev \cite{Lebedev3} in the proof.

\bigskip

{\bf Proof of Theorem \ref{t:main_intr}.}
One can suppose that $N$ is a sufficiently large prime number. Let
\begin{equation}\label{f:Q_def}
    Q = \exp ( (\log N)^\alpha (\log \log N)^{\beta}) \,,
\end{equation}
\begin{equation}\label{f:Theta_def}
    \Theta(x) = o ( (\log x)^c (\log \log x)^d ) \,,
\end{equation}
and
\begin{equation}\label{f:la_def}
    \nu = \Theta (Q) = o( (\log N)^{\alpha c} (\log \log N)^{\beta c + d} ) \,.
\end{equation}
Here $\alpha>1,\beta, c>0, d$ are real numbers which we will choose later.
Let $S$ be a subset of $\Z_N$.
Put $\Phi_{S,n} (x) = e^{2\pi i n \_phi(x)} \chi_S (x)$, $n= 1,\dots, Q$.
For any positive integer $k$ Lemma~\ref{T_k_est} gives us
\begin{equation}\label{tmp:26.11.2013_1}
    \oT_k ( \Phi_{S,n} ) \ge \frac{|S|^{2k}}{N \nu^{2k-2}} \,.
\end{equation}
Note that if $\| u \| > 2^{-1} Q^{-1/2}$ then the estimate $|\sin x| \ge 2x/\pi$, $x\in [0,\pi/2]$ implies
$$
    |\sum_{n=1}^Q e^{2\pi i nu}| \le 2^{-1} Q^{1/2} < Q^{1/2} \,.
$$
Using the last bound and summing inequality (\ref{tmp:26.11.2013_1}) over $n \in \Z_N$ it is easy to see that
the number of tuples $(x_1,\dots,x_{2k}) \in S^{2k}$
such that
$$
    x_1 + \dots + x_k = x_{k+1} + \dots + x_{2k} \,, \quad
$$
\begin{equation}\label{f:tuples}
    \| \_phi^* (x_1) + \dots + \_phi^* (x_k) - \_phi^* (x_{k+1}) - \dots - \_phi^* (x_{2k}) \| \le 2^{-1} Q^{-1/2}
\end{equation}
is at least $\frac{|S|^{2k}}{2N \nu^{2k-1}}$,
provided by the following holds
\begin{equation}\label{cond:Q_simple}
    Q \ge 4N^{6k} \ge 4 |S|^{-2} N^2 \nu^{4k-4} \,.
\end{equation}
In other words, for any subset $S\subseteq \Z_N$, we have $\oT^{\_phi,2^{-1} Q^{-1/2}}_k (S) \ge \frac{|S|^{2k}}{2N \nu^{2k-1}}$.

Now consider a maximal subset $\L \subseteq \Z_N$ belonging to the family
$\L^{\_phi,Q^{-1/2}} (2k,s)$, where the parameter $s$ depends on the constants $\alpha,\beta,c,d$ and
will be chosen later.
By the maximality of $\L$ we have that for any $x\in \Z_N$ there are elements $\l_1, \dots, \l_w \in \L$, $w\le 4k$
and integers $s_*, s_1,\dots,s_w$, $|s_*| \le s$, $s_* \neq 0$, $|s_j|\le s$, $j=1,\dots,w$,
$\sum_{j=1}^{|\L|} |s_j| + |s_*| \le 4k$ such that
\begin{equation}\label{tmp:05.12.2013_1}
    \| s_* \_phi^* (x) - \sum_{j=1}^w s_j \_phi^* (\l_j) \| \le Q^{-1/2} \,.
\end{equation}
Let us estimate the cardinality of $\L$.
Take $k=2 + [\log N]$.
Apply  Lemma \ref{l:Lambda_k_phi} to the set $\L$ and suppose that the first term in the maximum of formula
 (\ref{f:T_k_up_estimate_phi}) dominates.
Combining the estimate of the lemma and inequality $\oT^{\_phi,2^{-1} Q^{-1/2}}_k (\L) \ge \frac{|\L|^{2k}}{2N \nu^{2k-1}}$, we obtain
\begin{equation}\label{tmp:05.12.2013_2}
    |\L|
        \ll \nu^2 k \ll \nu^2 \log N \,.
\end{equation}
If the second term in the maximum of formula (\ref{f:T_k_up_estimate_phi}) dominates and
bound (\ref{tmp:05.12.2013_2}) does not hold then
$$
    (\log N)^{\frac{s}{s-4}} \gg |\L| \gg \nu^2 \log N
        \gg
            (\log N)^{1+2\alpha c} (\log \log N)^{2\beta c + 2d} \,.
$$
Choosing the parameter $s=s(\alpha,\beta,c,d)$ in such a way that the last inequality fails, we get that estimate (\ref{tmp:05.12.2013_2}) takes place
anyway.

Using Dirichlet theorem for elements of the set $\{ \_phi^* (x) \}_{x\in \L}$, we find a positive integer $q$ such that
$\| q\_phi^* (\l) \| \le 1/(4Nks!)$, $\l \in \L$ and $q\le (4Nks!)^{|\L|}$.
If
\begin{equation}\label{tmp:06.12.2013_1}
    |\L| \log N \ll \nu^2 \log^2 N \ll \log Q
\end{equation}
then choosing a sufficiently small constant  under $(\log x)^c (\log \log x)^d$ in $o$ from
formula (\ref{f:Theta_def}), we can assume that $q\le Q$ and, simultaneously,
\begin{equation}\label{tmp:17.12.2013_1}
    s! q Q^{-1/2} \le 1/2N \,.
\end{equation}
In terms of the constants  $\alpha,\beta,c,d$ inequality (\ref{tmp:06.12.2013_1}) can be rewritten as
\begin{equation}\label{tmp:06.12.2013_2}
    2\alpha c + 2 \le \alpha \quad \mbox{ and } \quad 2\beta c + 2d \le \beta \,.
\end{equation}
Hence if (\ref{tmp:06.12.2013_2}) takes place then inequalities (\ref{tmp:05.12.2013_1}) and (\ref{tmp:17.12.2013_1})
imply that for any  $x\in \Z_N$ the following holds
$$
    \| q s! \_phi^* (x) \|
        \le
            s! |s_* |^{-1} \| q s_* \_phi^* (x) \|
                \le
$$
$$
                \le
                    s! \| \sum_{j=1}^w s_j q \_phi^* (\l_j) \|+ s! q Q^{-1/2}
                        \le
                    1/2N + s! q Q^{-1/2} \le 1/N \,.
$$
In the last estimate we have used the fact that for any $0\neq s_*$, $|s_*| \le s$ the number $s! |s_* |^{-1}$
is a positive integer.

Let $\eta = 1/(64 \Theta^2 (Q))$.
Suppose that $\eta \ge (\log N)^{-1/4}(\log\log N)^{1/2}$.
Using Lemma \ref{l:Lebedev_extract} with parameter $Q=q$, we obtain
\begin{equation}\label{tmp:06.12.2013_3}
    \Theta (Q) \gg (\log N)^{\alpha c} (\log \log N)^{\beta c + d} \gg (\log N)^{1/10} (\log \log N)^{-3/10} \,.
\end{equation}
 Now choosing $\alpha = 2.2$, $c=1/22$, $\beta=-3/5$, $d=-3/11$, we have in view of (\ref{cond:main_22_1}) (or  (\ref{f:Theta_def})) a contradiction with (\ref{tmp:06.12.2013_3}).
  Simultaneously, we satisfy  condition (\ref{tmp:06.12.2013_2}).
Also it is easy to check that the choice implies that the condition $\eta \ge (\log N)^{-1/4}(\log\log N)^{1/2}$ takes place.
Finally, returning to (\ref{cond:Q_simple}) and recalling $k=2+[\log N]$ we see  that the inequality holds.
This completes the proof.
$\hfill\Box$

\bigskip

{\bf Remark.}
If one does not want to care about double logarithms in the main result then the using of Lemma \ref{l:Lambda_k_phi} can be avoided.
To do this one should consider {\it non--special } tuples $(x_1,\dots,x_{2k})$ that is the set of tuples satisfying (\ref{f:tuples})
with the condition any such tuple contains an element  $x_j$  appears once
(clearly, the number of the tuples do not have the property does not exceed $(Ck)^k |S|^k$, $C>0$ is an absolute constant).
Further, successfully remove such tuples to the moment when its number became less than $O(\Theta^2 (Q) k)$.
Then it is easy to see that for all $x\in \Z_N$ the quantity  $\_phi^* (x)$ can be expressed (with the accuracy $Q^{-1/2}$) as a combination of numbers $\_phi^* (x_j)$ with coefficients $\pm 1$,
where $x_j$ belong to the set of
non--special
tuples.
The remaining arguments repeat the proof of Theorem \ref{t:main_intr}.






{}
\bigskip

\noindent{S.V.~Konyagin\\
Steklov Mathematical Institute,\\
ul. Gubkina, 8, Moscow, Russia, 119991}
\\
and
\\
MSU,\\
Leninskie Gory, Moscow, Russia, 119992\\
{\tt konyagin@mi.ras.ru}

\bigskip

\noindent{I.D.~Shkredov\\
Steklov Mathematical Institute,\\
ul. Gubkina, 8, Moscow, Russia, 119991}
\\
and
\\
IITP RAS,  \\
Bolshoy Karetny per. 19, Moscow, Russia, 127994\\
{\tt ilya.shkredov@gmail.com}


\begin{thebibliography}{99}

\bibitem{BH}
{\sc A.~Beurling, H.~Helson. } Fourier--Stieltjes transforms with bounded powers // Math. Scand., {\bf 1}, 120--126, 1953.

\bibitem{GK}
{\sc B.J.~Green, S.V.~Konyagin. } On the Littlewood problem modulo a prime // Canad. J. Math. 2009. Vol.~61~(1), P.~141–-164.

\bibitem{Zyg}
{\sc A.~Zygmund. } Triginometric series / V.~2, CUP, 2002.

\bibitem{Kahane}
{\sc J.-P.~Kahane. } Transform\'{e}es de Fourier des fonctions sommables // Proceedings of the Int.
Congr. Math., in Inst. Mittag--Leffler, Djursholm, Sweden, 1963, pp. 114--131.

\bibitem{Kahane_book}
{\sc J.-P.~Kahane. } Absolutely convergent Fourier series / Mir, M., 1976.

\bibitem{Kon}
{\sc S.V.~Konyagin. } On a problem of Littlewood // Izvestiya of Russian Academy of Sciences, {\bf 45}:2, (1981), 243--265. English transl.: Math. USSR-Izv., 18:2 (1982), 205--225. 

\bibitem{Lebedev1}
{\sc V.V.~Lebedev. } Quantitative estimates in Beurling--Helson type theorems // Mat. Sbornik,
{\bf 201}:12, 103--130, 2010.

\bibitem{Lebedev2}
{\sc V.V.~Lebedev. } Estimates in Beurling--Helson type theorems: multidimensional case // Mat. Notes, {\bf 90}:3, 394--407, 2011.

\bibitem{Lebedev3}
{\sc V.V.~Lebedev. } Absolutely convergent Fourier series. An improvement of Beurling--Helson theorem. // Funct. Anal. Appl., {\bf 46}:2, 52--65, 2012.

\bibitem{MGPS}
{\sc O.C.~McGehee, L.~Pigno, B.~Smith. } Hardy's inequality and the $L^1$ norm of exponential sums //
Annals of Math., {\bf 113}, 613--618, 1981.

\bibitem{Rudin_book}
{\sc W.~Rudin. } Fourier analysis on groups /  Wiley 1990 (reprint of the 1962 original).

\bibitem{Sanders}
{\sc T.~Sanders.} The Littlewood--Gowers problem // J. Anal. Math., {\bf 101} 123--162, 2007.

\bibitem{Shkr_large_exp}
{\sc I.D.~Shkredov. } On Sets of Large Exponential Sums // Izvestiya of Russian Academy of Sciences, {\bf 72}:1, 161--182, 2008.


\end{thebibliography}
\end{document}